\RequirePackage{ifpdf}
\ifpdf 
\documentclass[pdftex]{sigma}
\else
\documentclass{sigma}
\fi

\newcommand{\R}{\mathbb{R}}
\renewcommand{\S}{{\boldsymbol{S}}}
\newcommand{\N}{\mathbb{N}}
\newcommand{\D}{{\mathcal D}}
\newcommand{\nn}{{{\boldsymbol{n}},{\boldsymbol{p}}}}
\newcommand{\E}{{\mathcal E}}
\renewcommand{\H}{{\mathcal H}}
\newcommand{\alphabold}{{\boldsymbol{\alpha}}}
\newcommand{\p}{{\boldsymbol{p}}}
\newcommand{\x}{{\boldsymbol{x}}}
\newcommand{\y}{{\boldsymbol{y}}}
\newcommand{\z}{{\boldsymbol{z}}}
\newcommand{\n}{{\boldsymbol{n}}}
\newcommand{\0}{{\boldsymbol{0}}}
\numberwithin{equation}{section}

\begin{document}
\allowdisplaybreaks

\renewcommand{\PaperNumber}{071}

\renewcommand{\thefootnote}{$\star$}

\FirstPageHeading

\ShortArticleName{Generalized Ellipsoidal and Sphero-Conal Harmonics}

\ArticleName{Generalized Ellipsoidal and Sphero-Conal Harmonics\footnote{This paper is a contribution 
to the Vadim Kuznetsov Memorial Issue ``Integrable Systems and Related Topics''.
The full collection is available at 
\href{http://www.emis.de/journals/SIGMA/kuznetsov.html}{http://www.emis.de/journals/SIGMA/kuznetsov.html}}}

\Author{Hans VOLKMER}
\AuthorNameForHeading{H. Volkmer}

\Address{Department of Mathematical Sciences, University of Wisconsin-Milwaukee,\\
P.\,O.~Box 413, Milwaukee, WI 53201 USA}
\Email{\href{mailto:volkmer@uwm.edu}{volkmer@uwm.edu}}
\URLaddress{\href{http://www.uwm.edu/~volkmer/}{http://www.uwm.edu/\~{}volkmer/}}

\ArticleDates{Received August 25, 2006, in f\/inal form October 20,
2006; Published online October 24, 2006}

\Abstract{Classical ellipsoidal and sphero-conal harmonics are polynomial solutions of
the Laplace equation that can be expressed in terms of Lam\'e polynomials.
Generalized ellipsoidal and sphero-conal harmonics are polynomial solutions of
the more general Dunkl equation that can be expressed in terms of Stieltjes polynomials.
Niven's formula
connec\-ting ellipsoidal and sphero-conal harmonics is generalized. Moreover, generalized ellipsoidal
harmonics are applied to solve
the Dirichlet problem for Dunkl's equation on ellipsoids.}


\Classification{33C50; 35C10}

\section{Introduction}
The theory of ellipsoidal and sphero-conal harmonics is a beautiful achievement of
classical mathematics.
It is not by accident that the well-known treatise ``A Course in Modern Analysis'' by
Whittaker and Watson \cite{WW} culminates in the f\/inal chapter ``Ellipsoidal Harmonics and
Lam\'e's Equation''.

An ellipsoidal harmonic is a polynomial $u(x_0,x_1,\dots,x_k)$ in $k+1$ variables
$x_0,x_1,\dots,x_k$
which satisf\/ies the Laplace equation
\begin{equation}\label{1:Laplace}
\Delta u:=\sum_{j=0}^k \frac{\partial^2 u}{\partial x_j^2}=0
\end{equation}
and assumes the product form
\begin{equation}\label{1:LameProduct1}
u(x_0,x_1,\dots,x_k)=E(t_0)E(t_1)\cdots E(t_k)
\end{equation}
in ellipsoidal coordinates $(t_0,t_1,\dots,t_k)$ with
$E$ denoting a Lam\'e quasi-polynomial.

A sphero-conal harmonic is
a polynomial $u(x_0,x_1,\dots,x_k)$
which satisf\/ies the Laplace equation~\eqref{1:Laplace}
and assumes the product form
\begin{equation}\label{1:LameProduct2}
u(x_0,x_1,\dots,x_k)=r^mE(s_1)E(s_2)\cdots E(s_k)
\end{equation}
in sphero-conal coordinates $(r,s_1,s_2,\dots,s_k)$.
Again, $E$ is a Lam\'e quasi-polynomial.

In most of the literature, for example, in the books by Hobson \cite{Hob} and
Whittaker and Watson~\cite{WW},
ellipsoidal and sphero-conal harmonics
are treated as polynomials of only three variables. The generalization to
any number of variables is straight-forward.
Since we plan to work in arbitrary dimension we will employ
ellipsoidal and sphero-conal coordinates in algebraic form; see Sections 3 and 4.
In $\R^3$ we may uniformize these coordinates by using Jacobian elliptic functions.
Sphero-conal harmonics are special instances of spherical harmonics. Spherical
harmonics in any dimension can be found in books by Hochstadt \cite{Ho} and
M\"uller \cite{Mu}.

We will extend the theory of ellipsoidal and sphero-conal harmonics
by replacing the Laplace equation by the equation
\begin{equation}\label{1:Dunkl}
    \Delta_\alphabold u:=\sum_{j=0}^k \D_j^2 u=0,\qquad \alphabold=(\alpha_0,\alpha_1,\dots,\alpha_k)
\end{equation}
introduced by Dunkl \cite{Du}.
In \eqref{1:Dunkl} we use the generalized partial derivatives
\begin{equation}\label{1:partial}
 \D_j u(\x):=\frac{\partial}{\partial x_j} u(\x)+
\alpha_j\frac{u(\x)-u(\sigma_j \x)}{x_j},\qquad \x=(x_0,x_1,\dots,x_k),
\end{equation}
where
$\sigma_j$ is the ref\/lection at the $j$th coordinate plane:
\[ \sigma_j(x_0,x_1,\dots,x_k)=(x_0,x_1,\dots,x_{j-1},-x_j,x_{j+1},\dots,x_k).\]
Equation \eqref{1:Dunkl} contains real parameters $\alpha_0,\alpha_1,\dots,\alpha_k$.
If $\alpha_j=0$ for all~$j$ then the equation reduces to the Laplace
equation.

A {\it generalized ellipsoidal harmonic} is a polynomial $u(x_0,x_1,\dots,x_k)$
which satisf\/ies Dunkl's equation \eqref{1:Dunkl}
and assumes the product form \eqref{1:LameProduct1}
in ellipsoidal coordinates
but with $E$ now denoting a Stieltjes quasi-polynomial. Generalized ellipsoidal harmonics
will be treated in Section~3 while
Stieltjes quasi-polynomials are introduced in Section 2.

\looseness=1
 A {\it generalized sphero-conal harmonic} is
a polynomial $u(x_0,x_1,\dots,x_k)$ which satisf\/ies Dunkl's equation \eqref{1:Dunkl}
and assumes the product form
\eqref{1:LameProduct2}
in sphero-conal coordinates. Again, $E$ is a~Stieltjes quasi-polynomial.
Generalized sphero-conal harmonics will be considered in Section~4.

It is very pleasing to see Stieltjes quasi-polynomials
taking over the role of Lam\'e quasi-polynomials.
Stieltjes polynomials have been considered for a long time but they did not appear in the
context of separated solutions of the Laplace equation.
Therefore, our paper shows how Stieltjes polynomials become part of the theory of
``Special Functions''.

It is quite remarkable that all the known results for classical ellipsoidal and sphero-conal
harmonics carry over to their generalizations.
In Section~5 we generalize formulas due to Hobson \cite[Chapter 4]{Hob}.
In  Section~6, as a consequence, we prove a generalization of Niven's formula
\cite[Chapter 23]{WW} connecting
ellipsoidal and sphero-conal harmonics.
In Section~7 we apply generalized ellipsoidal harmonics in order to solve a Dirichlet problem
for \eqref{1:Dunkl} on ellipsoids. This generalizes the classical result that
ellipsoidal harmonics may be used to
f\/ind the harmonic function which has prescribed values on the
boundary of an ellipsoid.
Finally, we give some examples in Section~8.

We point out that parts of Section 4 overlap with the author's paper \cite{Vo}.
The contents of the present paper are also related to the book by
Dunkl and Xu \cite{DuXu}, and the papers by Liamba  \cite{L} and Xu \cite{Xu},
however, these works do not involve Stieltjes polynomials.
We also refer to papers by Kalnins and Miller \cite{KM1,KM2,KM3}. The paper \cite{KM3}
addresses Niven's formula from a dif\/ferent perspective.
Kutznetsov \cite{Ku} and Kutznetsov and Komarov \cite{KK} have also worked in related areas.
Kutznetsov jointly with Sleeman wrote the chapter on Heun functions for the
Digital Library of Mathematical Functions. Stieltjes polynomials appear in this chapter.

\section{Stieltjes quasi-polynomials}
We consider the Fuchsian dif\/ferential equation
\begin{equation}\label{2:Fuchs}
\prod_{j=0}^k (t-a_j)\left[ v''+\sum_{j=0}^k \frac{\alpha_j+\frac12}
{t-a_j} v'\right] +\left[-\frac12\sum_{j=0}^k \frac{p_j\alpha_j A_j}{t-a_j}
+\sum_{i=0}^{k-1}\lambda_i t^i\right]v =0
\end{equation}
for the function $v(t)$ where the prime denotes dif\/ferentiation with respect to~$t$.
This dif\/ferential equation contains four sets of real parameters:
\begin{eqnarray}\label{2:a}
&& a_0<a_1<\dots<a_k,\\
\label{2:alpha}
&& \alpha_0,\alpha_1,\dots,\alpha_k\in(-\tfrac12,\infty),\\
\label{2:epsilon}
&& p_0,p_1,\dots,p_k\in\{0,1\},\\
\label{2:lambda}
&& \lambda_0,\lambda_1,\dots\lambda_{k-1}\in\R,
\end{eqnarray}
and $A_j$ is an abbreviation:
\begin{equation}\label{2:A}
 A_j:=\prod_{i=0 \atop i\ne j}^k (a_j-a_i).
\end{equation}
Usually, the f\/irst three sets of parameters are given while the $\lambda$'s
play the role of eigenvalue parameters.

Equation \eqref{2:Fuchs} has regular singularities at inf\/inity and at each $a_j$, $j=0,1,\dots,k$.
The exponents at $a_j$ are $\nu_j=\frac{p_j}2$ and $\mu_j=\frac{1-p_j}{2}-\alpha_j$.
If $\nu_{k+1}$, $\mu_{k+1}$ denote the exponents at inf\/inity then
\begin{equation}\label{2:index1}
 \lambda_{k-1}=\nu_{k+1}\mu_{k+1}
\end{equation}
 and
\begin{equation}\label{2:index2}
\sum_{j=0}^{k+1}(\nu_j+\mu_j)=k.
\end{equation}
The accessory parameters $\lambda_0,\lambda_1,\dots,\lambda_{k-2}$ are unrelated to
the exponents.

The following result def\/ines {\it  Stieltjes quasi-polynomials} $E_\nn$.\
Let parameter sets \eqref{2:a}  and \eqref{2:alpha} be given.
For every multi-index $\n=(n_1,n_2,\dots,n_k)\in\N^k$ of nonnegative integers
and $\p=(p_0,p_1,\dots,p_k)\in\{0,1\}^k$ there exist uniquely determined values of the parameters
$\lambda_0,\dots,\lambda_{k-1}$ such that \eqref{2:Fuchs} admits a solution
of the form
\begin{equation}\label{2:quasi}
 E_\nn(t)=\left(\prod_{j=0}^k |t-a_j|^{p_j/2}\right) \tilde E_\nn(t),\qquad t\in\R,
\end{equation}
where $\tilde E_\nn$ is a polynomial with exactly
$n_j$ zeros in the open interval $(a_{j-1},a_j)$ for
each $j=1,\dots,k$. The polynomial $\tilde E_\nn(t)$ is uniquely determined
up to a constant factor and has the degree $|\n|=n_1+\dots+n_k$.
We normalize $\tilde E_\nn$ so that its leading coef\/f\/icient is unity.
Then we may write $E_\nn$ in the form
\begin{equation}\label{2:quasi2}
E_\nn(t)=\prod_{j=0}^k |t-a_j|^{p_j/2}
\prod_{\ell=1}^{|\n|} (t-\theta_\ell),
\end{equation}
where
\[ \theta_1<\theta_2<\dots<\theta_{|\n|}. \]
Then $\theta_1,\dots\theta_{n_1}$ lie in $(a_0,a_1)$, $\theta_{n_1+1},\dots,\theta_{n_1+n_2}$
lie in $(a_1,a_2)$ and so on.

If $\p=\0$ then $E_{\n,\0}$ is a polynomial introduced by Stieltjes
\cite{St} whose work was inf\/luenced by Heine \cite[Part III]{He}.
A proof of existence and uniqueness of the polynomials $E_{\n,\mathbf 0}$
can be found in Szeg\"o \cite[Section 6.8]{Sz}.
For general $\p$ a computation shows that $\tilde E_\nn$
is the Stieltjes polynomial~$E_{\n,\0}$ with $\alpha_j$ replaced by $\alpha_j+p_j$.
Therefore, the proof of existence and uniqueness in the
general case can be reduced to the special case $\p=\0$.

The value of $\lambda_{k-1}$ associated with $E_\nn$ can be computed.
One of the exponents at inf\/inity must be
\[ \nu_{k+1}=|\n|+\tfrac12|\p|, \qquad |\p|:=\sum_{j=0}^k p_j. \]
Using \eqref{2:index1}, \eqref{2:index2}, we obtain
\begin{equation}\label{2:lambdak}
\lambda_{k-1}=-\tfrac12m\left(\tfrac12m+|\alphabold|+\tfrac{k-1}{2}\right),\qquad m:=2|\n|+|\p|,\qquad
|\alphabold|:=\sum_{j=0}^k \alpha_j .
\end{equation}
No formulas are known for the corresponding values of the accessory parameters.

If $\alpha_j=0$ for all $j$ then Stieltjes quasi-polynomials reduce to
Lam\'e quasi-polynomials in arbitrary dimension. If we work in
$\R^3$ there are eight possible choices of the parameters
\eqref{2:epsilon} giving us the familiar eight types of classical
Lam\'e quasi-polynomials; see Arscott \cite{A}.

\section{Generalized ellipsoidal harmonics}

We say that a function $u:\R^{k+1}\to\R$ has parity
$\p=(p_0,\dots,p_k)\in\{0,1\}^{k+1}$ if
\[ u(\x)-u(\sigma_j\x)=2p_j u(\x) \qquad\text{for \ \ $j=0,1,\dots,k$}.\]

Equation \eqref{1:Dunkl} can be written in the form
\begin{equation}\label{3:Dunkl}
\Delta u(\x) +\sum_{j=0}^k \frac{2\alpha_j}{x_j}\frac{\partial}{\partial x_j} u(\x)
-\sum_{j=0}^k \frac{\alpha_j}{x_j^2}(u(\x)-u(\sigma_j \x))=0 .
\end{equation}
If $u$ has parity $\p$
then \eqref{3:Dunkl} becomes the partial dif\/ferential equation
\begin{equation}\label{3:Dunkl2}
\Delta u(\x) +\sum_{j=0}^k \frac{2\alpha_j}{x_j}\frac{\partial}{\partial x_j}u(\x)
-\sum_{j=0}^k \frac{2p_j \alpha_j}{x_j^2}u(\x)=0 .
\end{equation}

In order to introduce ellipsoidal coordinates,
f\/ix the parameters \eqref{2:a}.
For every $(x_0,\dots,x_k)$ in the positive cone of $\R^{k+1}$
\begin{equation}\label{3:cone}
x_0>0,\dots,x_k>0,
\end{equation}
its ellipsoidal coordinates $t_0,t_1,\dots,t_k$ lie in the intervals
\begin{equation}\label{3:cube}
 a_k<t_0<\infty,\qquad a_{i-1}<t_i<a_i,\qquad i=1,\dots,k,
\end{equation}
and satisfy
\begin{equation}\label{3:ellipsoidal1}
\sum_{j=0}^k \frac{x_j^2}{t_i-a_j} =1\qquad\text{for \ \ $i=0,1,\dots,k$}.
\end{equation}
Conversely, for given $t_i$ in the intervals \eqref{3:cube},
we have
\begin{equation}\label{3:ellipsoidal2}
    x_j^2=\frac{\prod\limits_{i=0}^k (t_i-a_j)}{\prod\limits_{i=0\atop i\ne j}^k
    (a_i-a_j)} .
\end{equation}
These coordinates provide a bijective mapping between the positive
cone \eqref{3:cone} and the cube~\eqref{3:cube}.

We now transform the partial dif\/ferential equation \eqref{3:Dunkl2} for functions
$u(\x)$ def\/ined on the cone~\eqref{3:cone}
to ellipsoidal coordinates, and then we apply the method of separation of variables.
We obtain $k+1$ times the Fuchsian equation~\eqref{2:Fuchs}
coupled by the separation constants $\lambda_0,\lambda_1,\dots,\lambda_{k-1}$.
We do not carry out the details of these known calculations.
A good reference is Schmidt and Wolf~\cite{SW}.
Therefore, if $v_j(t_j)$, $j=0,1\dots,k$,  are solutions of \eqref{2:Fuchs}
with $t_j$ ranging in the intervals \eqref{3:cube}
then the function
\begin{equation}\label{3:separation}
 u(x_0,\dots,x_k)=v_0(t_0)\cdots v_k(t_k)
\end{equation}
satisf\/ies \eqref{3:Dunkl2}.
Of course,
the values of the parameter sets \eqref{2:a}--\eqref{2:lambda}
must be the same in each equation \eqref{2:Fuchs}.

As a special case choose $v_j$ as the Stieltjes quasi-polynomial
$E_\nn$ for each~$j$. Then we know that
\begin{equation}\label{3:product}
  F_\nn(x_0,x_1,\dots,x_k):=E_\nn(t_0)E_\nn(t_1)\cdots E_\nn(t_k)
\end{equation}
solves \eqref{3:Dunkl2}.
This function $F_\nn$ is our generalized ellipsoidal harmonic.

\begin{theorem}\label{3:t1}
The generalized ellipsoidal harmonic $F_\nn$ is a polynomial in $x_0,x_1,\dots,x_k$ which satisfies
Dunkl's equation \eqref{1:Dunkl}. It is of total degree $2|\n|+|\p|$
and has parity~$\p$.
\end{theorem}
\begin{proof}
If $t_0,\dots,t_k$ denote ellipsoidal coordinates of
$x_0,\dots,x_k$, then
\begin{equation}\label{3:ellipsoidal3}
   \prod_{j=0}^k (t_j-\theta)=\left(\prod_{i=0}^k (a_i-\theta)\right)
\left(1-\sum_{j=0}^k \frac{x_j^2}{\theta-a_j}\right)
\end{equation}
for every $\theta$ dif\/ferent from each $a_j$. In fact,
both sides of \eqref{3:ellipsoidal3} are polynomials in $\theta$
of degree $k+1$ with leading coef\/f\/icient $(-1)^{k+1}$. Moreover, both sides of the
equation vanish at $\theta=t_0,t_1,\dots,t_k$ by def\/inition \eqref{3:ellipsoidal1}.
So equation \eqref{3:ellipsoidal3} follows.

By \eqref{2:quasi2}, \eqref{3:ellipsoidal2}, \eqref{3:product} and \eqref{3:ellipsoidal3}, the function
$F_\nn$ can be written as
\begin{equation}\label{3:F}
   F_\nn(\x)=c_\nn \x^\p
\prod_{\ell=1}^{|\n|}\left(\sum_{j=0}^k \frac{x_j^2}{\theta_\ell-a_j}-1\right),
\end{equation}
where
\[ \x^\p:=x_0^{p_0}\cdots x_k^{p_k}\qquad\text{for} 
  \quad \x=(x_0,\dots,x_k),\quad \p=(p_0,\dots,p_k) ,\]
and $c_\nn$ is the constant
\begin{equation}\label{3:c}
 c_\nn:=(-1)^{|\n|}\left(\prod_{j=0}^k |A_j|^{p_j/2}\right)
\left(\prod_{\ell=1}^{|\n|}\prod_{i=0}^k(a_i-\theta_\ell)\right)
\end{equation}
with $A_j$ according to \eqref{2:A}. This shows that $F_\nn$
is a polynomial of total degree
$2|\n|+|\p|$.
We know that $F_\nn$ solves \eqref{3:Dunkl2} on the cone \eqref{3:cone} and since it has parity $\p$
it solves~\eqref{1:Dunkl} on~$\R^{k+1}$.
\end{proof}

\section{Generalized sphero-conal harmonics}
In order to introduce sphero-conal coordinates,
f\/ix the parameters \eqref{2:a}. Let
$(x_0,x_1,\dots,x_k)$ be in the positive cone \eqref{3:cone} of
$\R^{k+1}$. Its sphero-conal coordinates $r, s_1,\dots,s_k$ are
determined in the intervals
\begin{equation}\label{4:cube}
 r>0,\qquad a_{i-1}<s_i<a_i,\qquad i=1,\dots,k
\end{equation}
by the equations
\begin{equation}\label{4:sphero1}
 r^2=\sum_{j=0}^k x_j^2
\end{equation}
and
\begin{equation}\label{4:sphero2}
 \sum_{j=0}^k \frac{x_j^2}{s_i-a_j}=0\qquad\text{for \ \ $i=1,\dots,k$}.
\end{equation}
This def\/ines a bijective map from the positive cone in $\R^{k+1}$ to
the set of points $(r,s_1,\dots,s_k)$ satisfying~\eqref{4:cube}.
The inverse map is given by
\begin{equation}\label{4:inverse}
  x_j^2=r^2 \frac{\prod\limits_{i=1}^k(s_i-a_j)}{\prod\limits_{i=0 \atop i\neq j}^k (a_i-a_j)}.
\end{equation}

We now transform the partial dif\/ferential equation \eqref{3:Dunkl2} for functions
$u(x_0,x_1,\dots,x_k)$ def\/ined on the cone \eqref{3:cone}
to sphero-conal coordinates and then we apply the method of separation of variables
\cite{SW}.
For the variable $r$ we obtain the Euler equation
\begin{equation}\label{4:Euler}
 v_0''+\frac{k+2|\alpha|}{r} v_0'+
\frac{4\lambda_{k-1}}{r^2}v_0=0
\end{equation}
while for the variables $s_1,s_2,\dots,s_k$  we obtain the Fuchsian equation \eqref{2:Fuchs}.
More precisely, if $\lambda_0,\ldots,\lambda_{k-1}$
are any given numbers (separation constants),
if $v_0(r)$, $r>0$, solves \eqref{4:Euler}
and $v_i(s_i)$, $a_{i-1}<s_i<a_i$, solve \eqref{2:Fuchs}
for each $i=1,\dots,k$, then
\[ u(x_0,x_1,\dots,x_k)=v_0(r)v_1(s_1)v_2(s_2)\cdots v_k(s_k) \]
solves \eqref{3:Dunkl2}.

Let $E_\nn$ be a Stieltjes quasi-polynomial.
It follows from \eqref{2:lambdak} that
$v_0(r)=r^m$ is a solution of \eqref{4:Euler}, where $m:=2|\n|+|\p|$.
Therefore,
\begin{equation}\label{4:product}
 G_\nn(x_0,x_1,\dots,x_k):=r^m E_\nn(s_1)E_\nn(s_2)\cdots  E_\nn(s_k)
\end{equation}
is a solution of \eqref{3:Dunkl2}.
This function $G_\nn$ is our generalized sphero-conal harmonic.

\begin{theorem}\label{4:t1}
The generalized sphero-conal harmonic $G_\nn$
is a polynomial in $x_0,x_1,\dots,x_k$, it is homogeneous of degree $2|\n|+|\p|$,
it has parity $\p$ and it solves Dunkl's equation \eqref{1:Dunkl}.
\end{theorem}
\begin{proof}
Let $(x_0,\dots,x_k)$ be
a point with $x_j>0$ for all $j$, and let $(r,s_1,\dots,s_k)$
denote its corresponding sphero-conal coordinates. We claim that
\begin{equation}\label{4:eq1}
 r^2(s_1-\theta)\dots(s_k-\theta)=\left(\prod_{i=0}^k (a_i-\theta)\right)
\sum_{j=0}^k \frac{x_j^2}{a_j-\theta}
\end{equation}
for all $\theta$ which are dif\/ferent from each $a_j$.
Both sides of \eqref{4:eq1} are polynomials in $\theta$ of degree~$k$ with leading
coef\/f\/icient $(-1)^kr^2$. Moreover, both sides vanish at $\theta=s_1,\dots,s_k$
because of def\/inition \eqref{4:sphero2}.
Equation \eqref{4:eq1} is established.

We write the Stieltjes quasi-polynomial $E_\nn$ in the form
\eqref{2:quasi2}.
Using \eqref{4:inverse}, \eqref{4:product} and \eqref{4:eq1}, we obtain
\begin{equation}\label{4:G}
   G_\nn(\x)=c_\nn \x^\p
     \prod_{\ell=1}^{|\n|} \sum_{j=0}^k \frac{x_j^2}{\theta_\ell-a_j} ,
\end{equation}
where $c_\nn$ is given by \eqref{3:c}.
This shows that $G_\nn(\x)$ is a polynomial in $x_0,x_1,\dots,x_k$, it is
homogeneous of degree
$2|\n|+|\p|$, and it has parity $\p$.
We know that $G_\nn$ solves \eqref{3:Dunkl2} on the cone \eqref{3:cone} and since it has
parity $\p$
it solves
\eqref{1:Dunkl} on $\R^{k+1}$.
\end{proof}

A {\it generalized spherical harmonic} is a homogeneous polynomial $u$ in the variables $x_0,x_1,\dots,x_k$
which solves Dunkl's equation
\eqref{1:Dunkl}.
For a given set of parameters  \eqref{2:alpha} we let
$\H_m$ denote the f\/inite dimensional linear space of all generalized spherical harmonics of degree $m$.
If $\alpha_j=0$ for each $j$ then we obtain the classical spherical harmonics.

On the $k$-dimensional unit sphere $\S^k$ we introduce the inner product
\begin{equation}\label{4:inner}
 \langle f,g\rangle_w :=\int_{\S^k}  w(\x)f(\x)g(\x)\,dS(\x),
\end{equation}
and norm
\begin{equation}\label{4:norm}
\|f\|_w:=\langle f,f\rangle_w^{1/2},
\end{equation}
where the weight function $w$ is def\/ined by
\begin{equation}\label{4:weight}
w(x_0,x_1,\dots,x_k):=|x_0|^{2\alpha_0}|x_1|^{2\alpha_1}\cdots|x_k|^{2\alpha_k}.
\end{equation}
The surface measure on the sphere is normalized so that $\int_{\S^k} dS(\x)$ equals the
surface area of the sphere $\S^k$.
The condition $\alpha_j>-\frac12$ ensures that $\langle f,g\rangle_w$ is well-def\/ined if
$f$ and $g$ are continuous on $\S^k$.

\begin{theorem}\label{4:t2}
Let $m\in\N$.
The system of all generalized sphero-conal harmonics $G_\nn$ of degree~$m$
forms an orthogonal basis for $\H_m$ with respect to the inner product
\eqref{4:inner}.
\end{theorem}
\begin{proof}
We consider the system
of all sphero-conal harmonics $G_\nn$, where $\nn$ satisfy
$m=2|\n|+|\p|$.
By Theorem \ref{4:t1}, $G_\nn$ belongs to $\H_m$.
The dimension of the linear space of generalized spherical harmonics of degree $m$ which
have parity $\p$
is
\begin{equation}\label{4:dim}
 \binom{\frac12m-\frac12|\p|+k-1}{k-1}
\end{equation}
if $m-|\p|$ is a nonnegative even integer and zero otherwise.
This can be proved as in Hochstadt \cite[p.~170]{Ho} or it follows from
Dunkl \cite[Proposition 2.6]{Du2}
where a basis of $\H_m$ in terms of Jacobi polynomials is constructed.
The dimension \eqref{4:dim} agrees with the number of multi-indices $\n=(n_1,n_2,\dots,n_k)$ for which
$m=2|\n|+|\p|$.
We conclude that the number of pairs $\nn$ with $m=2|\n|+|\p|$ agrees
with the dimension of $\H_m$.
Therefore, in order to complete the proof of the theorem, we have to show that
$G_\nn$ is orthogonal
to $G_{\n',\p'}$ provided $(\n,\p)\ne
(\n',\p')$.
If $\p\ne\p'$ this is clear because the weight function \eqref{4:weight}
is  an even function.
If $\p=\0$ and $\n\ne \n'$ then orthogonality was shown in \cite[Theorem 3.3]{Vo}.
The proof of orthogonality in the remaining cases is analogous  and is omitted.
\end{proof}

Extending the method of proof of \cite[Theorem 3.3]{Vo} we also
establish the following theorem.
\begin{theorem}\label{4:t3}
The system of all generalized sphero-conal harmonics $G_\nn$, $\n\in\N^k$,
$\p\in\{0,1\}^{k+1}$, when properly normalized, forms
an orthonormal basis of $L^2_w(\S^k)$.
\end{theorem}

No explicit formula is known for the norm of $G_\nn$ in $L^2_w(\S^k)$.
However, the norm of a~polynomial can be computed using
the formula
\begin{equation}\label{4:gamma}
\int_{\S^k} |x_0|^{2\beta_0-1}|x_1|^{2\beta_1-1}\cdots|x_k|^{2\beta_k-1}\,dS(\x)=
\frac{2\prod\limits_{j=0}^k\Gamma(\beta_j)}{\Gamma(\beta_0+\beta_1+\dots+\beta_k)}
\end{equation}
which holds whenever $\beta_j>0$, $j=0,1,\dots,k$.

In general, ellipsoidal harmonics are not homogeneous polynomials so they are not spherical harmonics.
However, they are related to spherical harmonics in the
following way.

\begin{theorem}\label{4:t4}
Let $m\in\N$. The system of all generalized ellipsoidal harmonics $F_\nn$ of total degree
$2|\n|+|\p|$ at most $m$ is a basis for the direct sum
\begin{equation}\label{4:direct}
 \H_0 \oplus \H_1\oplus\dots \oplus \H_m .
\end{equation}
\end{theorem}
\begin{proof}
The Dunkl operator $\Delta_\alphabold$ maps a homogeneous polynomial
of degree $q$ to a homogeneous polynomial of degree $q-2$.
Therefore, if we write a generalized  ellipsoidal harmonic as a sum of homogeneous polynomials,
then these homogeneous polynomials also satisfy \eqref{1:Dunkl}.
Hence every generalized ellipsoidal harmonic of total degree at most $m$ lies in the direct sum
\eqref{4:direct}.
By comparing \eqref{3:F} with \eqref{4:G}, we see that
\begin{equation}\label{4:compare}
 F_\nn(\x)=G_\nn(\x)+\text{terms of lower degree.}
\end{equation}
By Theorem \ref{4:t2}, the system of all $G_\nn$ with $2|\n|+|\p|\le m$
is a basis for the direct sum \eqref{4:direct}.
The statement of the theorem follows.
\end{proof}

Of course, the spaces $\H_m$ of generalized
spherical harmonics and the direct sum \eqref{4:direct} depend on the parameters
$\alpha_0,\alpha_1,\dots,\alpha_k$. However, it is easy to show that the
set of functions on $\S^k$ which are restrictions of function in \eqref{4:direct}
is independent of these parameters. In fact, this set consists of all functions that are
restrictions of polynomials of total degree at most $m$ to~$\S^k$.

\section{Hobson's formulas}
In this section we generalize some formulas given by Hobson \cite[p.\ 124]{Hob}.
These formulas will be applied in the next section to obtain a generalization of Niven's formula.
\begin{lemma}\label{5:hobson1}
Let $\D$ be the operator given by
\[ \D f(x):= f'(x)+\alpha\frac{f(x)-f(-x)}{x}, \]
where $\alpha$ is a constant.
Then, for $A(z):=z^\ell$, $\ell,m\in\N$, we have
\begin{equation}\label{5:hob1}
 \D^m x^{2\ell}=\sum_{j=0}^m 2^{m-2j} A^{(m-j)}\big(x^2\big)\frac{1}{j!}\D^{2j} x^m .
\end{equation}
\end{lemma}
\begin{proof}
If $m>2\ell$, then both sides of \eqref{5:hob1} are zero. So we
assume that $m\le 2\ell$. We f\/irst consider the case that $m=2n$ is even.
The left-hand side of \eqref{5:hob1} is equal to
\begin{equation}\label{5:hob2}
2^m \frac{\ell!}{(\ell-n)!}(-1)^n\big(\tfrac12-\ell-\alpha\big)_n x^{2\ell-m} .
\end{equation}
The right-hand side of \eqref{5:hob1} is equal to
\begin{equation}\label{5:hob3}
\sum_{j=0}^n 2^{m-2j}\frac{\ell!}{(\ell-m+j)!} \frac1{j!}
2^{2j} \frac{n!}{(n-j)!}(-1)^j\big(\tfrac12-\alpha-n\big)_j\, x^{2\ell-m} .
\end{equation}
After some simplif\/ications, equality of \eqref{5:hob2} and \eqref{5:hob3}
follows from the Chu-Vandermonde sum
\[ \sum_{j=0}^n \frac{(a)_j(b)_{n-j}}{j!(n-j)!}=\frac{(a+b)_n}{n!} \]
applied to $a=\frac12-\alpha-n$, $b=n-\ell$.
This completes the proof of \eqref{5:hob1} if $m$ is even.
The similar proof for odd $m$ is omitted.
\end{proof}

Clearly, in \eqref{5:hob1} it would be enough to let $j$ run from $0$ to
$\lfloor \frac m2 \rfloor$. Similar remarks apply to other formulas in this section.

In the following lemma, $\D_j$ is according to \eqref{1:partial} and
$\partial_j$ is the usual partial derivative with respect to $x_j$.

\begin{lemma}\label{5:hobson2}
Let $m_0,m_1,\dots,m_k\in \N$, and let $A:(0,\infty)^{k+1}\to\R$
be $m:=m_0+m_1+\dots+m_k$ times differentiable.
Then, for $x_0,\dots,x_k\neq 0$,
\begin{gather}
\D_0^{m_0}\cdots\D_k^{m_k} [A(x_0^2,\dots,x_k^2)]
=\sum_{j_0=0}^{m_0}\dots\sum_{j_k=0}^{m_k} 2^{m-2(j_0+\dots+j_k)}\nonumber \\
\qquad{}\times(\partial_0^{m_0-j_0}\cdots\partial_k^{m_k-j_k}A)(x_0^2,\dots,x_k^2)
\frac{\D_0^{2j_0}\cdots\D_k^{2j_k}}{j_0!\cdots j_k!}
x_0^{m_0}\cdots x_k^{m_k} .\label{5:hob4}
\end{gather}
Warning: On the left-hand side of this formula the operators $\D_j$ are
applied to the function $f(x_0,\dots,x_k):=A(x_0^2,\dots,x_k^2)$,
whereas on the right-hand side the partial derivatives $\partial_j$
are applied directly to $A$.
\end{lemma}
\begin{proof}
Let $B$ be the Taylor polynomial of $A$ of order $m$ at a given point
$(z_0,\dots,z_k)$ with $z_j>0$. Let $x_j:=z_j^{1/2}$.
Then \eqref{5:hob4} is true with $A-B$ in place of $A$ at the point
$x_0,\dots,x_k$ (both sides of the equation are zero.)
Therefore, it is suf\/f\/icient to prove \eqref{5:hob4} for polynomials $A$,
and so for monomials
\[ A(z_0,\dots,z_k)=z_0^{\ell_0} \cdots z_k^{\ell_k} .\]
In this case, we obtain \eqref{5:hob4} by applying Lemma \ref{5:hobson1}
to each function $z_j^{\ell_j}$ and multiplying.
\end{proof}

If $f(x_0,\dots,x_k)$ is a polynomial, we will use the operator
$f(\D_0,\dots,\D_k)$. It is well-def\/ined because the operators $\D_j$ commute.
We use
\[ r:=(x_0^2+\dots+x_k^2)^{1/2} .\]

\begin{theorem}\label{5:hobson3}
Let $f_m(x_0,\dots,x_k)$ be a homogeneous polynomial of degree~$m$,
and let $B:(0,\infty)\to \R$ be $m$ times differentiable.
Then, for all nonzero $(x_0,x_1,\dots,x_k)$,
\begin{equation}\label{5:hob5}
f_m(\D_0,\dots,\D_k)[B(r^2)]=
\sum_{j=0}^m 2^{m-2j} B^{(m-j)}(r^2)\frac1{j!} \Delta_\alphabold^j
f_m(x_0,\dots,x_k) .
\end{equation}
\end{theorem}
\begin{proof}
It is suf\/f\/icient to prove \eqref{5:hob5} for monomials
\[ f_m(x_0,\dots,x_k)=x_0^{m_0}\cdots x_k^{m_k},\qquad m=m_0+\cdots+m_k .\]
In this case \eqref{5:hob5} follows from Lemma
\ref{5:hobson2} with $A(z_0,\dots,z_k)=B(z_0+\cdots+z_k)$
by using
\[ \partial_0^{m_0-j_0}\cdots \partial_k^{m_k-j_k} A(z_0,\dots,z_k)
=B^{(m-j)}(z_0+\dots+z_k) \]
with $j=j_0+\dots+j_k$ and the multinomial formula
\begin{gather*}
 \Delta_\alphabold^j=(\D_0^2+\dots+\D_k^2)^j=
\sum_{j_0+\dots+j_k=j} \binom j {j_0 \cdots j_k} \D_0^{2j_0}\cdots\D_k^{2j_k} .\tag*{\qed}
\end{gather*}\renewcommand{\qed}{}
\end{proof}

When we apply Theorem \ref{5:hobson3} to $B(z):=z^\gamma$ with $\gamma\in\R$,
we obtain the following corollary.

\begin{corollary}\label{5:hobson4}
Let $f_m(x_0,\dots,x_k)$ be a homogeneous polynomial of degree~$m$.
Then, for $\gamma\in \R$,
\begin{equation}\label{5:hob6}
f_m(\D_0,\dots,\D_k) [r^{2\gamma}]
= \sum_{j=0}^m 2^{m-2j}(-1)^{m-j}(-\gamma)_{m-j} r^{2(\gamma-m+j)}
\frac1{j!}\Delta_\alphabold^j f_m(x_0,\dots,x_k) .
\end{equation}
\end{corollary}

\begin{corollary}\label{5:hobson5}
Let $f_m$ be a generalized spherical harmonic of degree $m$.
Then, for $\gamma\in\R$,
\begin{equation}\label{5:hob8}
r^{2(m-\gamma)} f_m(\D_0,\dots,\D_k) [r^{2\gamma}]
= 2^m (-1)^m(-\gamma)_m f_m(x_0,\dots,x_k) .
\end{equation}
\end{corollary}

If we set
\begin{equation}\label{5:gamma}
\gamma=\frac{1-k}2-|\alphabold|,
\end{equation}
then a simple calculation shows that
$\Delta_\alphabold r^{2\gamma}=0$.
So $r^{2\gamma}$ plays the role of a fundamental solution
of $\Delta_\alphabold u=0$ generalizing the
solution $1/r$ of the Laplace equation $\Delta u=0$ in $\R^3$
with $\alpha=(0,0,0)$.
Note that the number $\gamma$ def\/ined by \eqref{5:gamma} is always less than $1$.
It can be zero (for example for the Laplacian in the plane). In this case,
$\ln r$ plays the role of a fundamental solution.
The fundamental solution $r^{2\gamma}$ and an associated formula producing harmonic polynomials
appeared in Xu \cite{Xu2}.

\begin{corollary}\label{5:hobson6}
If $f_m(x_0,\dots,x_k)$ is a homogeneous polynomial of degree $m$ and
$\gamma$ is defined by \eqref{5:gamma}, then the right-hand side of
equation~\eqref{5:hob6} is a generalized spherical harmonic of degree~$m$.
\end{corollary}
\begin{proof}
Since $\Delta_\alphabold r^{2\gamma}=0$, this follows by applying
$\Delta_\alphabold$ to both sides of \eqref{5:hob6}.
\end{proof}

\section{Niven's formula}

In this section we prove a generalization of Niven's formula
expressing ellipsoidal harmonics in terms of sphero-conal harmonics.
We follow the method of Hobson \cite[p.\ 483]{Hob}.

Let $E=E_\nn$ be a Stieltjes quasi-polynomial, and let $F=F_\nn$, $G=G_\nn$ be the corresponding
ellipsoidal and sphero-conal harmonics written in the forms \eqref{3:F} and \eqref{4:G},
respectively.
It will be convenient to introduce the auxiliary polynomial
\[ H(\x):=c_\nn \x^\p
\prod_{\ell=1}^{|\n|} \left(\sum_{j=0}^k \frac{x_j^2}{\theta_\ell-a_j}-
\sum_{j=0}^k \frac{x_j^2}{t-a_j}\right), \]
where $t$ is a f\/ixed number greater than $a_k$.
We def\/ine positive constants $d_j$ by $d_j^2=t-a_j$ for $j=0,\dots,k$.
Let $\gamma$ be the constant def\/ined by \eqref{5:gamma}. We assume that
$\gamma\neq 0$.
The identity
\[ \sum_{j=0}^k \frac{x_j^2}{\theta-a_j} -\sum_{j=0}^k
\frac{x_j^2}{t-a_j}= (t-\theta)
\sum_{j=0}^k \frac{x_j^2}{(\theta-a_j)(t-a_j) }, \]
implies
\begin{equation}\label{6:niven3}
 H(d_0x_0,\dots,d_kx_k)=E(t)G(x_0,\dots,x_k).
\end{equation}
By Corollary \ref{5:hobson5}, we have
\begin{equation}\label{6:niven4}
2^m (-1)^m (-\gamma)_m G(x_0,\dots,x_k)=r^{2(m-\gamma)}G(\D_0,\dots,\D_k)
[r^{2\gamma}] ,
\end{equation}
where $m:=2|\n|+|\p|$.
Since
\[ r^2+(t-\theta)\sum_{j=0}^k \frac{x_j^2}{\theta-a_j}
=\sum_{j=0}^k \frac{(t-a_j)x_j^2}{\theta-a_j},\]
we conclude that
\[ G(d_0x_0,\dots,d_kx_k)=E(t)G(x_0,\dots,x_k)+r^2P(x_0,\dots,x_k), \]
where $P$ is a polynomial.
It follows that
\[ G(d_0\D_0,\dots,d_k\D_k)=
E(t) G(\D_0,\dots,\D_k)+P(\D_0,\dots,\D_k)\Delta_\alphabold .\]
Using that $\Delta_\alphabold r^{2\gamma}=0$, we obtain
\begin{equation}\label{6:niven5}
G(d_0\D_0,\dots,d_k\D_k)\big[r^{2\gamma}\big]=
E(t) G(\D_0,\dots,\D_k)\big[r^{2\gamma}\big].
\end{equation}

We now combine equations \eqref{6:niven3}, \eqref{6:niven4}, \eqref{6:niven5}
and obtain
\begin{equation}\label{6:niven6}
2^m(-1)^m(-\gamma)_m
H(d_0x_0,\dots,d_kx_k)=r^{2(m-\gamma)}G(d_0\D_0,\dots,d_k\D_k)\big[r^{2\gamma}\big].
\end{equation}
Now \eqref{6:niven6} and Corollary \ref{5:hobson4} yield
\[ H(d_0x_0,\dots,d_kx_k)=\sum_{i=0}^m \frac{r^{2i}}{2^{2i}i!
(\gamma-m+1)_i}
\Delta_\alphabold^i [G(d_0x_0,\dots,d_kx_k) ] .\]
We replace the variables $x_j$ by $y_j/d_j$ and rename $y_j$ as
$x_j$ again. This gives
\begin{equation}\label{6:niven7}
 H(\x)=\sum_{i=0}^m \frac{R^{2i}}{2^{2i}i!(\gamma-m+1)_i}
\big(d_0^2\D_0^2+\dots+d_k^2\D_k^2\big)^i G(\x) ,
\end{equation}
where
\[ R^2=\sum_{j=0}^k \frac{x_j^2}{d_j^2}=\sum_{j=0}^k\frac{x_j^2}{t-a_j} .\]
Since $G$ satisf\/ies $\Delta_\alphabold G=0$, we see
that the right-hand side of \eqref{6:niven7}
does not change if we replace $d_j^2\D_j^2$ by
$-a_j\D_j^2$.
If positive numbers $x_0,\dots,x_k$ are given, we can choose
$t$ as the ellipsoidal coordinate $t=t_0$. Then $R=1$ and
$H(\x)=F(\x)$.
We have proved the following theorem.

\begin{theorem}\label{6:niven}
Let $F_\nn,G_\nn$ be the generalized ellipsoidal and sphero-conal harmonics
of degree $m=2|\n|+|\p|$
defined by \eqref{3:F}, \eqref{4:G}, respectively.
Assume that $\gamma$ def\/ined by \eqref{5:gamma} is nonzero.
Then
\begin{equation}\label{6:niven8}
F_\nn(\x)=\sum_{i=0}^m \frac{(-1)^i}{2^{2i}i!(\gamma+1-m)_i}
\big(a_0\D_0^2+\dots+a_k\D_k^2\big)^i G_\nn(\x) .
\end{equation}
\end{theorem}

In the classical case $k=2$ and $\alpha_0=\alpha_1=\alpha_2=0$ this is Niven's
formula; see \cite[p.\ 489]{Hob}.

\section{A Dirichlet problem for ellipsoids}
In this section we apply generalized ellipsoidal harmonics to solve the Dirichlet
boundary value problem for the Dunkl equation on ellipsoids.

We consider the solid ellipsoid
\[ \E:=\left\{\x\in\R^{k+1} : \sum_{j=0}^k
\frac{x_j^2}{b_j^2}<1 \right\},
\]
with semi-axes $b_0>b_1>\dots>b_k>0$.
Let $\partial\E$ be the boundary of $\E$.
Given a function $f:\partial\E\to\R$ we want to f\/ind a solution $u$
of Dunkl's equation \eqref{1:Dunkl} on $\E$ that assumes the given boundary values
$f$ on $\partial\E$ in the sense explained below.

It will be convenient to parameterize $\partial\E$ by the unit sphere $\S^k$
employing the map
\begin{equation}\label{7:T}
T: \S^k \to \partial \E
\end{equation}
def\/ined by
\[ T(y_0,y_1,\dots,y_k):=(b_0y_0,b_1y_1,\dots,b_ky_k).
\]
We suppose that the given boundary value function $f:\partial \E\to\R$
has the property that the function $f\circ T$ is in
$L^2_w(\S^k)$,
where the weight function $w$ is def\/ined in \eqref{4:weight}.
A solution of the Dirichlet boundary value problem for the Dunkl equation
with given boundary value function~$f$ is a function  $u\in C^2(\E)$
which satisf\/ies \eqref{1:Dunkl} in $\E$
and assumes the boundary value $f$ in the
following sense. For suf\/f\/iciently small $\delta>0$, form the
confocal ellipsoids
\begin{equation}\label{7:confocal}
\left\{\x\in\R^{k+1}: \sum_{j=0}^k \frac{x_j^2}
{b_j^2-\delta}=1 \right\},
\end{equation}
and let $T_\delta$ be def\/ined as $T$ but with respect to the ellipsoid \eqref{7:confocal}
in place
of $\partial \E$. Then we require that
\begin{equation}\label{7:bc}
 u\circ T_\delta\to f\circ T\qquad
\text{in \ \ $L^2_w(\S^k)$ \ \ as \ \ $\delta\to 0$}.
\end{equation}

We now show how to construct a solution of this Dirichlet problem.
We choose any real number $\omega$ (we can take $\omega=0$ if we wish),
and def\/ine numbers $a_0<a_1<\dots<a_k$ by
\[ a_j:=\omega-b_j^2 .\]
Corresponding to these numbers $a_j$
we introduce sphero-conal coordinates $(r,s_1,\dots,s_k)$ for cartesian coordinates
$\y=(y_0,y_1,\dots,y_k)$ and ellipsoidal
coordinates $(t_0,t_1,\dots,t_k)$ for cartesian coordinates $\x=(x_0,x_1,\dots,x_k)$.
Note that if $\x=T\y$ and $r=1$ then $s_j=t_j$ for each $j=1,2,\dots,k$.

Since the function $f\circ T$ lies in $L_w^2(\S^k)$,
we can expand $f\circ T$
in the orthonormal basis of Theorem \ref{4:t3}:
\begin{equation}\label{7:Fourier}
 f\circ T=\sum_\nn f_\nn e_\nn G_{n,p},
\end{equation}
where the factors $e_\nn$ are determined by
\[ e_\nn\|G_\nn\|_w=1,   \]
and
\begin{equation}\label{7:d}
f_\nn:=\langle f\circ T,e_\nn G_\nn\rangle_w.
\end{equation}
Then
\begin{equation}\label{7:d2}
\sum_\nn |f_\nn|^2=\|f\circ T\|_w^2<\infty .
\end{equation}
The expansion \eqref{7:Fourier} converges in $L^2_w(\S^k)$.
We are going to prove that
\begin{equation}\label{7:sol}
u(\x):=\sum_{\nn} \frac{f_\nn e_\nn}{E_\nn(\omega)}F_\nn(\x)
\end{equation}
is the desired solution of our Dirichlet problem.

\begin{theorem}\label{7:t1}
The function $u$ defined by \eqref{7:sol} is infinitely many times differentiable and solves
Dunkl's equation \eqref{1:Dunkl} on the open ellipsoid $\E$, and it assumes the given
boundary value $f$ in the sense of \eqref{7:bc}.
\end{theorem}
\begin{proof}
We f\/irst show that the inf\/inite series in \eqref{7:sol} converges.
We know from \cite[Lemma~2.2]{L} that there is a sequence $K_m$ of polynomial growth such that
\begin{equation}\label{7:ineq1}
|e_\nn G_\nn(\y)|\le  K_m\qquad\text{for \ \ $m=2|\n|+|\p|$, \ \ $\y\in \S^k$.}
\end{equation}
For given $t\in(a_k,\omega)$ we consider the solid ellipsoid
\[ \E_t:=\left\{\x\in\R^{k+1} : \sum_{j=0}^k\frac{x_j^2}{t-a_j}\le 1\right\}\]
which is a subset of $\E$.
By comparing \eqref{3:product} and \eqref{4:product}, we get from \eqref{7:ineq1}
\begin{equation}\label{7:ineq2}
|e_\nn F_\nn(\x)|\le   E_\nn(t)K_m\qquad\text{for \ \ $\x\in \partial\E_t$.}
\end{equation}
The Stieltjes quasi-polynomial $E_\nn$ has degree $m/2$ and all of its zeros lie
in the interval $[a_0,a_k]$. Hence we have the inequality
\begin{equation}\label{7:ineq3}
0<E_\nn(t)\le E_\nn(\omega)\left(\frac{t-a_0}{\omega-a_0}\right)^{m/2}
\qquad \text{for \ \ $a_k<t\le\omega$}.
\end{equation}
Now we obtain from \eqref{7:ineq2}, \eqref{7:ineq3}
\begin{equation}\label{7:ineq4}
 \left|e_\nn F_\nn(\x)\right|
\le \left(\frac{t-a_0}{\omega-a_0}\right)^{m/2}E_\nn(\omega)K_m\qquad\text{for \ \ $x\in \E_t$}.
\end{equation}
Since the set of numbers $f_\nn$ is bounded by \eqref{7:d2},
$K_m$ grows only polynomially with $m$ and
$\big(\frac{t-a_0}{\omega-a_0}\big)^{m/2}$ goes to $0$ exponentially as $m\to\infty$, we see that
the series in \eqref{7:sol} converges uniformly in $\E_t$ and thus in
every compact subset of $\E$.

The next step is to show that $u$ solves equation \eqref{1:Dunkl}.
This follows if we can justify interchanging the operator $\Delta_\alphabold$ with the sum in
\eqref{7:sol}.
Consider the series
\begin{equation}\label{7:sol2}
 \sum_{\nn}\frac{f_\nn e_\nn}{E_\nn(\omega)}\frac{\partial}{\partial x_j}
F_\nn(x_0,x_1,\dots,x_k)
\end{equation}
that we obtain from \eqref{7:sol} by dif\/ferentiating each term with
respect to $x_j$.
In order to show uniform convergence of this series on $\E_t$
we need a bound for the partial
derivatives of $F_\nn$. We obtain such a bound from the following
result due to Kellog \cite{Ke}. If $P(y_0,y_1,\dots,y_k)$ is a polynomial of total degree
$N$ then
\begin{equation}\label{7:Kellog}
\max\{\|{\rm grad}\, P(\y)\|: \|y\|\le 1\}\le N^2 \max\{|P(\y)|: \|y\|\le1 \} ,
\end{equation}
where $\|\cdot\|$ denotes euclidian norm in $\R^{k+1}$.
If we use a mapping like \eqref{7:T} to transform the ellipsoid to the unit ball
we f\/ind that
\begin{equation}\label{7:ineq5}
\left|\frac{\partial}{\partial x_j}F_{\nn}(\x)\right|\le \frac{m^2}{(t-a_k)^{1/2}}
\max\{| F_\nn(\z)|: \z\in\E_t \}\qquad\text{for \ \ $\x\in\E_t$.}
\end{equation}
Using this estimate
we show as before that \eqref{7:sol2}
converges uniformly on $\E_t$.
In a similar way we argue for the second term in the generalized partial derivative
\eqref{1:partial}.
Since we can repeat the procedure we see
that $u$ is inf\/initely many times dif\/ferentiable on $\E$ and it solves
Dunkl's equation.

It remains to show that $u$ satisf\/ies the boundary condition \eqref{7:bc}.
For $y\in \S^k$, we f\/ind
\[ f\circ T(\y)-f\circ T_\delta(\y)=\sum_{\nn} f_\nn e_\nn
\left(1-\frac{E_\nn(\omega-\delta)}{E_\nn(\omega)}
\right) G_\nn(\y). \]
Hence we obtain
\[ \|f\circ T-f\circ T_\delta\|_w^2=\sum_{\nn}f_\nn^2\left(1-\frac{E_\nn(\omega-\delta)}
{E_\nn(\omega)}
\right)^2 \]
so \eqref{7:bc} follows easily.
\end{proof}

\section{Examples}
Formulas in this paper have been checked with the software {\it Maple}
for some Stieltjes polynomials represented in explicit form.
For example, take $k=2$,
\[ a_0=0,\qquad a_1=3,\qquad a_2=5,\quad
\alpha_0=\frac{229}{54},\qquad \alpha_1=\frac{71}{54},\qquad \alpha_2=\frac{25}{6},\]
and
\[ n_1=2,\qquad n_2=1,\qquad p_0=p_1=p_2=0.\]
Then the corresponding Stieltjes polynomial $E_{\n,\p}$
is given by
\[ E_{\n,\p}(t)=(t-1)(t-2)(t-4). \]
Indeed, $E_\nn$ satisf\/ies equation \eqref{2:Fuchs} with
\[ \lambda_0=\frac{1120}{9},\qquad \lambda_1=-\frac{119}{3} ,\]
and it has two zeros between $a_0$ and $a_1$, and one zero between $a_1$ and $a_2$.

The simplest way to compute such examples is to use the fact
that the zeros $\theta_\ell$ of $E_{\n,\0}$ are characterized by
the system of equations
\begin{equation}\label{8:system}
\sum_{q=1 \atop q\ne \ell}^{|\n|} \frac{2}{\theta_\ell-\theta_q} +\sum_{j=0}^k
\frac{\alpha_j+\frac12}{\theta_\ell-a_j} =0,\qquad\ell=1,2,\dots,|\n|;
\end{equation}
see \cite[(6.81.5)]{Sz}.

The corresponding ellipsoidal and sphero-conal harmonics are
\begin{gather*}
 F_\nn=-192\left( x_0^2-\tfrac12x_1^2-\tfrac14x_2^2-1 \right)
   \left(\tfrac12x_0^2-x_1^2-\tfrac13x_2^2-1 \right)
  \left(\tfrac14x_0^2+x_1^2-x_2^2-1 \right),\\
G_\nn=-192\left( x_0^2-\tfrac12x_1^2-\tfrac14x_2^2 \right)
   \left(\tfrac12x_0^2-x_1^2-\tfrac13x_2^2 \right)
  \left(\tfrac14x_0^2+x_1^2-x_2^2 \right).
\end{gather*}
One can check that these polynomials do satisfy equation \eqref{1:Dunkl}.
Also, applying formula \eqref{6:niven8} to $G_\nn$ we obtain $F_\nn$ as claimed.

We now take the same $a_j$ but replace the parameters $\alpha_j$ by
\[ \alpha_0=\frac{229}{54}-1=\frac{175}{54},\qquad \alpha_1=\frac{71}{54}-1=\frac{17}{54},\qquad
\alpha_3=\frac{25}{6} .\]
Moreover, let
\[ n_1=2,\qquad n_2=1,\qquad p_0=1,\qquad p_1=1,\qquad p_2=0 .\]
Then the Stieltjes quasi-polynomial is
\[ E_\nn(t)=\sqrt{|t|}\sqrt{|t-3|}(t-1)(t-2)(t-4) .\]
It satisf\/ies equation \eqref{2:Fuchs} with
\[ \lambda_0=\frac{2855}{18},\qquad \lambda_1=-\frac{440}{9} .\]
The corresponding ellipsoidal and sphero-conal harmonics are as before
but with
$-192$ replaced by $-192\sqrt{15}\sqrt{6}$ and the extra factor $x_1x_2$ added.
Again it can be checked that these polynomials satisfy equation \eqref{1:Dunkl},
and formula \eqref{6:niven8} holds.

\subsection*{Acknowledgements}
The author thanks W. Miller Jr.\ and two anonymous referees for helpful comments.

\LastPageEnding
\end{document}